\documentclass[11pt]{article}
\usepackage{amsfonts,amsmath,amsthm,amssymb}
\usepackage{graphicx}
\usepackage{booktabs,array}
\usepackage{needspace}
\usepackage{xcolor}
\usepackage{tikz}
\usepackage[hyperfootnotes=false,colorlinks=true,linkcolor=blue,citecolor=blue,urlcolor=black]{hyperref}

\allowdisplaybreaks
\numberwithin{equation}{section}
\newtheorem{theorem}{Theorem}[section]
\newtheorem{proposition}[theorem]{Proposition}
\newtheorem{lemma}[theorem]{Lemma}
\theoremstyle{definition}
\newtheorem{definition}[theorem]{Definition}
\newtheorem{example}[theorem]{Example}
\theoremstyle{remark}
\newtheorem{remark}[theorem]{Remark}
\newcommand{\Nzero}{\mathbb N_0}
\newcommand{\Sig}{\{1,2\}^{\Nzero}}
\newcommand{\cl}{\overline}
\newcommand{\Fol}{\operatorname{Fol}}
\newcommand{\suf}{\operatorname{suf}}
\newcommand{\Lang}{\mathcal L}
\newcommand{\pc}{\mathcal C}

\theoremstyle{plain}
\newtheorem{maintheorem}{Theorem}

\title{Subshifts of finite type and parameter density\\
for intermediate \(-\beta\)-transformations}
\author{Yiming Ding$^1$, Yun Sun$^{2,\dag}$\\[1ex]
\small\itshape $^1$School of Mathematics and Systems Science,\\
\small\itshape Wuhan University of Science and Technology, Wuhan 430081, China\\
\small\itshape $^2$ School of Mathematics and Statistics, Wuhan University of Technology,\\
\small\itshape Wuhan 430070, China\\[1ex]
\small E-mail: \texttt{dingym@wust.edu.cn}; \texttt{sunyun@whut.edu.cn}}
\date{}
\hypersetup{
  pdftitle={Subshifts of finite type and parameter density for intermediate negative beta-transformations},
  pdfauthor={Yiming Ding and Yun Sun},
  pdfkeywords={Intermediate negative beta-transformation; subshift of finite type; kneading invariants; parameter density}
}

\begin{document}
\maketitle
\begingroup
\renewcommand{\thefootnote}{}
\footnotetext{\(\dagger\) Corresponding author.}
\endgroup

\begin{center}
\begin{minipage}{120mm}
{\small\textbf{Abstract.}
We study the subshifts \(\Omega_{-\beta,\alpha}\) associated with the intermediate \(-\beta\)-transformations
\(x\mapsto-\beta x+\alpha\pmod1\), where \((\beta,\alpha)\in\Delta_-:=\{(b,a):1<b<2,\ b-1<a<1\}\). We prove that $\Omega_{-\beta,\alpha}$ is of
finite type if and only if both endpoints of \([0,1]\) reach the discontinuity \(c_{\beta,\alpha}=\alpha/\beta\).
We also prove that the set of parameters \((\beta,\alpha)\) for which
\(\Omega_{-\beta,\alpha}\) is a subshift of finite type is dense in \(\Delta_-\).
The approximating parameters can be chosen so that
both kneading invariants converge to those of the original map.
}
\end{minipage}
\end{center}
\vskip0.5cm

\section{Introduction}\label{sec:introduction}

The study of \(\beta\)-expansions goes back to the work of
R\'enyi~\cite{Renyi} and Parry~\cite{Parry}.   The associated symbolic space $\Omega_{\beta}$ is determined
by the expansion of the endpoint \(1\). A subshift which can be described by a finite set of forbidden
words is called a subshift of finite type. For convenience, we
write SFT. Parry proved that
the greedy \(\beta\)-shift is an SFT if and only if the
greedy expansion of \(1\) is finite. Such bases are called simple
\(\beta\)-numbers and form a dense subset of \((1,\infty)\). This result allows us
to approximate a general base by bases whose shifts are SFTs.
We refer the reader to Blanchard~\cite{Blanchard} and
Lind and Marcus~\cite{LM} for further background.

 For the intermediate $\beta$-transformations
\(x\mapsto\beta x+\alpha\pmod1\), where \((\beta,\alpha)\in\Delta:=\{(b,a):1<b<2,\ 0<a<2-b\}\), Li, Sahlsten and Samuel~\cite{LSS} proved
that the $\beta$-shift is an SFT if and only if both kneading
invariants are periodic. When considering density of SFT parameters,
Li et al.~\cite{LSSS} proved that these parameters are dense in $\Delta$. Sun, Li and Ding~\cite{SLD} further studied
density on fibers, including the case with one fixed
periodic kneading invariant.

For negative bases, Ito and Sadahiro~\cite{IS} introduced
\(-\beta\)-expansions. Here the symbolic comparisons use \emph{alternating lexicographic order},
since each branch reverses orientation. Liao and
Steiner~\cite{LS} studied the dynamical properties of
\((-\beta)\)-transformations.
Steiner~\cite{Steiner} gave necessary and sufficient
conditions for a sequence to be the \((-\beta)\)-expansion of \(1\)
for some \(\beta>1\).
Frougny and
Lai~\cite[Theorem~5.5]{FL} proved that the  \((-\beta)\)-shift
is an SFT if and only if the expansion of the   endpoint
\(-\beta/(\beta+1)\) is purely periodic.
Huang, Liao and Sun~\cite{HLS}
proved that the simple \(-\beta\)-numbers are dense in \((1,\infty)\).
Dajani and Ramawadh~\cite{DR} gave necessary and sufficient
conditions for soficity and finite type for a family of
$-\beta$-transformations.
Barnsley, Steiner and Vince~\cite{BSV} characterized
the critical itineraries of maps with constant slope and one discontinuity.
Other related results concern the specification property~\cite{HLY},
zeta functions~\cite{NNZ}, intrinsic ergodicity~\cite{SY}, and
measures of maximal entropy~\cite{NNBM}.

\begin{figure}[!htbp]
\centering
\resizebox{\linewidth}{!}{%
\begin{tikzpicture}[x=4.6cm,y=4.6cm,font=\small,
  mapline/.style={blue!70!black,line width=0.9pt},
  guide/.style={black,densely dashed,line width=0.4pt}]
  \pgfmathsetmacro{\comparisonbreak}{8/15}
  \foreach \panel in {0,1} {
    \begin{scope}[xshift=\panel*7.2cm]
      \draw[guide] (0,1)--(1,1)--(1,0);
      \draw[guide] (\comparisonbreak,0)--(\comparisonbreak,1);
      \draw[black,dashed,line width=0.45pt] (0,0)--(1,1);
      \draw[->,line width=0.5pt] (0,0)--(1.09,0)
        node[right] {$x$};
      \draw[->,line width=0.5pt] (0,0)--(0,1.10);
      \node[below left] at (0,0) {$0$};
      \draw (1,0.014)--(1,-0.014) node[below] {$1$};
      \draw (0.014,1)--(-0.014,1) node[left] {$1$};
    \end{scope}
  }
  \begin{scope}
    \node at (0.5,1.34) {(a) Positive slope};
    \node[above] at (0,1.10) {$S(x)$};
    \draw[guide] (0,0.2)--(\comparisonbreak,0.2);
    \draw[guide] (0,0.7)--(1,0.7);
    \draw (\comparisonbreak,0.014)--(\comparisonbreak,-0.014)
      node[below] {$p$};
    \draw (0.014,0.2)--(-0.014,0.2) node[left] {$a$};
    \draw (0.014,0.7)--(-0.014,0.7) node[left] {$b+a-1$};
    \draw[mapline] (0,0.2)--(\comparisonbreak,1);
    \draw[mapline] (\comparisonbreak,0)--(1,0.7);
    \fill[blue!70!black] (0,0.2) circle (1.5pt);
    \fill[blue!70!black] (1,0.7) circle (1.5pt);
    \draw[mapline,fill=white] (\comparisonbreak,1) circle (2.1pt);
    \draw[mapline,fill=white] (\comparisonbreak,0) circle (2.1pt);
  \end{scope}
  \begin{scope}[xshift=7.2cm]
    \node at (0.5,1.34) {(b) Negative slope};
    \node[above] at (0,1.10) {$T(x)$};
    \draw[guide] (0,0.8)--(\comparisonbreak,0.8);
    \draw[guide] (0,0.3)--(1,0.3);
    \draw (\comparisonbreak,0.014)--(\comparisonbreak,-0.014)
      node[below] {$c$};
    \draw (0.014,0.8)--(-0.014,0.8) node[left] {$\alpha$};
    \draw (0.014,0.3)--(-0.014,0.3)
      node[left] {$1+\alpha-\beta$};
    \draw[mapline] (0,0.8)--(\comparisonbreak,0);
    \draw[mapline] (\comparisonbreak,1)--(1,0.3);
    \fill[blue!70!black] (0,0.8) circle (1.5pt);
    \fill[blue!70!black] (1,0.3) circle (1.5pt);
    \draw[mapline,fill=white] (\comparisonbreak,0) circle (2.1pt);
    \draw[mapline,fill=white] (\comparisonbreak,1) circle (2.1pt);
  \end{scope}
\end{tikzpicture}}
\caption{Intermediate $\beta$- and $-\beta$-transformations with
$b=\beta=3/2$, $a=1/5$ and $\alpha=4/5$, so that
$p=c=8/15$. The branches have slope $+b$ in panel~(a) and
$-\beta$ in panel~(b).  }
\label{fig:positive-negative-maps}
\end{figure}
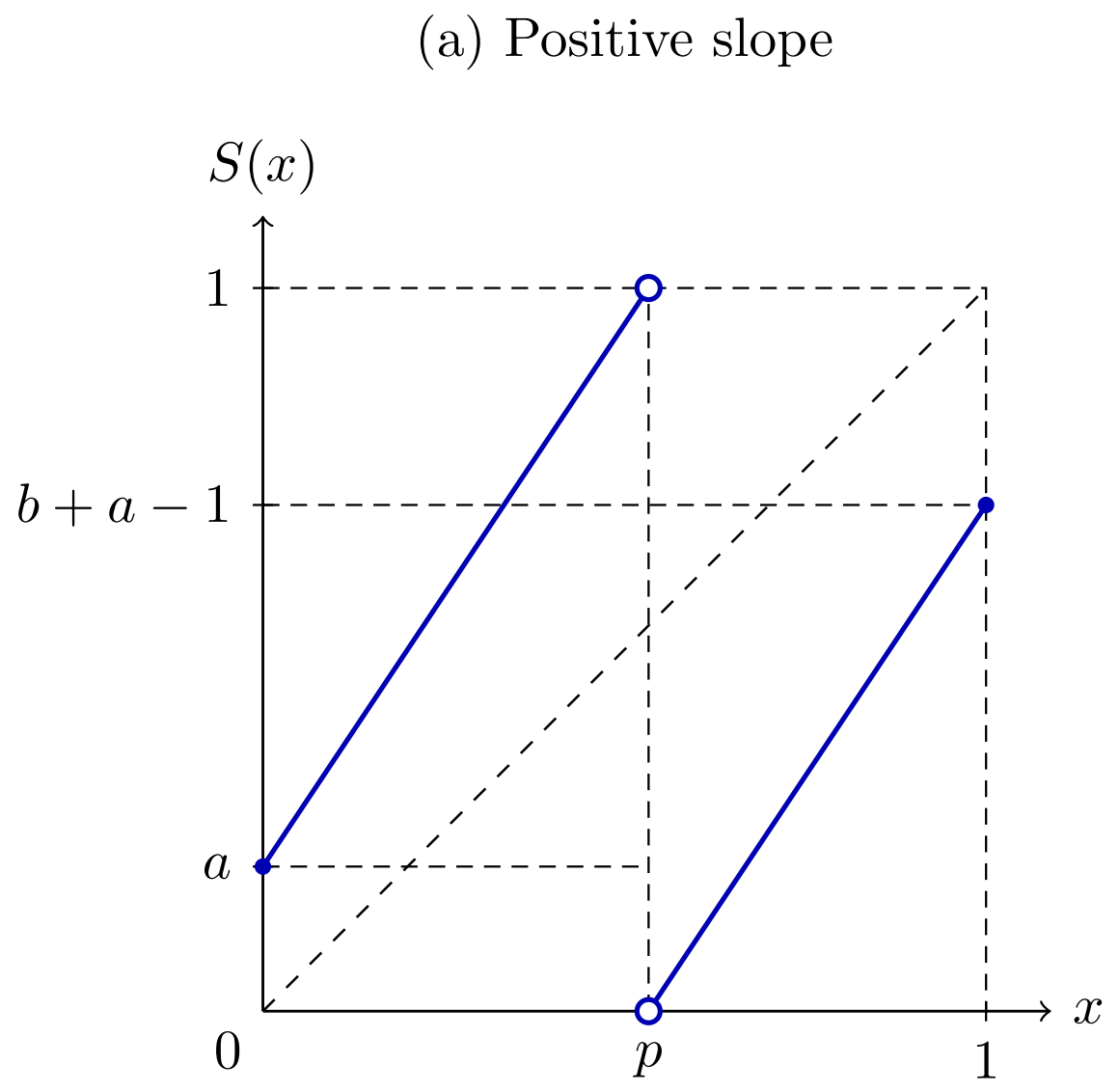

In this paper, we focus on the finite-type criterion and parameter
density for intermediate \(-\beta\)-transformations. Let
\(c=c_{\beta,\alpha}=\alpha/\beta\), and define
\begin{equation*}
 T(x)=T_{-\beta,\alpha}(x)=
 \begin{cases}
  \alpha-\beta x,&0\leq x<c,\\
  \alpha+1-\beta x,&c<x\leq1.
 \end{cases}
\end{equation*}
Then \(T(c-)=0\), \(T(c+)=1\), and both branches have slope
\(-\beta\). We restrict the
parameters to
\begin{equation*}
 \Delta_-:=\{(\beta,\alpha)\in\mathbb R^2:
                   1<\beta<2,\ \beta-1<\alpha<1\}.
\end{equation*}
This restriction forces the maps to have precisely two decreasing
branches and one discontinuity.   We use \(-\beta x+\alpha\pmod1\) for these two
branches, and specify the value at \(c\) when needed.
For the boundary case \(\alpha=1\), we have
\(c=1/\beta\) and set \(T(c)=1\). Then the restriction of \(T\)
to \((0,1]\) is affinely conjugate to the standard
Ito--Sadahiro \((-\beta)\)-transformation~\cite{IS,LS}.
For comparison, Figure~\ref{fig:positive-negative-maps} also shows
the positive-slope map \(S_{b,a}(x)=bx+a\pmod1\), with
discontinuity \(p=(1-a)/b\). Its branches preserve orientation,
while the branches of \(T\) reverse orientation.

We label the left and right branches by \(1\) and \(2\),
respectively, as in~\cite{DSS}. For a point \(x\) whose orbit avoids
\(c\), denote its kneading sequence by \(\tau(x)\). For
\(x\in[0,1]\), denote the one-sided limits through such points by
\(\tau(x-)\) and \(\tau(x+)\), taking only the inward limits
\(\tau(0+)\) and \(\tau(1-)\) at the endpoints. See Definition \ref{def:naturalshift} for details.
Throughout the paper, we call the
one-sided limits \(k_-=\tau(c-)\) and \(k_+=\tau(c+)\) the
\emph{kneading invariants}, following the convention of Lorenz
kneading theory~\cite{DSS,HS}. One-sided limit codings also appear
in Dajani and Ramawadh~\cite{DR}. The endpoint expansions
in~\cite{FL,HLS} and the critical itineraries in~\cite{BSV} use
fixed branch conventions at discontinuities. For negative slopes,
these codings need not coincide with the corresponding one-sided limits.
The \emph{intermediate \((-\beta)\)-shift} of \(T_{-\beta,\alpha}\) is defined by
\[
 \Omega_{-\beta,\alpha}
 =\{\tau(x-):0<x\leq1\}
 \cup\{\tau(x+):0\leq x<1\}
 \subset\{1,2\}^{\mathbb N_0}.
\]
This set is a subshift described by the kneading
invariants in Proposition~\ref{prop:kneadingdescription}.
  Notice that a return to \(c\) after an odd number of
iterates exchanges the side of approach, since each branch reverses
orientation. Hence a shift of one kneading invariant may equal the
other. Our first result gives necessary and sufficient conditions
for \(\Omega_{-\beta,\alpha}\) to be an SFT.

\begin{samepage}
\begin{maintheorem}\label{main:criterion}
For \((\beta,\alpha)\in\Delta_-\), the following are equivalent:
\begin{enumerate}
 \item \(\Omega_{-\beta,\alpha}\) is an SFT;
 \item both \(0\) and \(1\) reach \(c_{\beta,\alpha}\) in finite time;
 \item there exist \(p,q\geq1\) such that
 \(\sigma^pk_-\in\{k_-,k_+\}\) and
 \(\sigma^qk_+\in\{k_-,k_+\}\).
\end{enumerate}
\end{maintheorem}
\end{samepage}

The proof of Theorem~\ref{main:criterion} follows from
Propositions~\ref{prop:endpointreturn} and~\ref{prop:kneadingreturn}.
Notice that the kneading invariants need not both be purely
periodic. In Example~\ref{ex:mixed}, the shift is an SFT, but
\(k_-=112(2211)^\infty\) is not purely periodic.

To study density of SFT parameters, we introduce the notation
\[
 \mathcal F_-:=\{(\beta,\alpha)\in\Delta_-:
                                \Omega_{-\beta,\alpha}\text{ is an SFT}\}.
\]
It is natural to ask whether every parameter in \(\Delta_-\)
can be approximated by parameters in \(\mathcal F_-\).
The following theorem gives a positive answer.

\begin{maintheorem}\label{main:density}
The set \(\mathcal F_-\) is dense in \(\Delta_-\).
\end{maintheorem}

The proof of Theorem \ref{main:density} also gives convergence of the kneading invariants.
For every \((\beta,\alpha)\in\Delta_-\), we can choose
\((b_i,a_i)\in\mathcal F_-\) converging to \((\beta,\alpha)\)
such that \(k_-^{(i)}\to k_-\) and \(k_+^{(i)}\to k_+\),
where \((k_-^{(i)},k_+^{(i)})\) and \((k_-,k_+)\) are the
kneading invariants of \(T_{-b_i,a_i}\) and \(T_{-\beta,\alpha}\),
respectively.

The remainder of the paper is organized as follows.
In Section~\ref{sec:preliminaries}, we introduce
the symbolic notations, kneading invariants,
and recall some useful results. Section~\ref{sec:criterion}
 explains the role of return
parity and proves Theorem~\ref{main:criterion}.
In Section~\ref{sec:constructive}, we study endpoint
return curves and prove Theorem~\ref{main:density}.

\Needspace{7\baselineskip}
\section{Preliminaries}\label{sec:preliminaries}

\subsection{Subshifts}\label{subsec:subshifts}

We first recall some definitions and notations of
symbolic dynamics. We work with the space \(\Sig\),
where \(\Nzero=\{0,1,2,\ldots\}\), equipped with the usual metric:
\begin{equation*}
 d(\omega,\nu)=
 \begin{cases}
  0,&\omega=\nu,\\
  2^{-k},&k=\min\{j\geq0:\omega_j\ne\nu_j\},\quad\omega\ne\nu.
 \end{cases}
\end{equation*}
Note that the topology induced by this metric coincides with the
product topology, and \(\Sig\) is compact. The left-shift map
$
 \sigma(\omega_0\omega_1\omega_2\cdots)=\omega_1\omega_2\omega_3\cdots
$
is continuous.
Let \(\omega^{(n)}\in\Sig\) for \(n\geq1\).
Then \(\omega^{(n)}\to\omega\) if and only if,
for every \(L\geq1\), the first \(L\) symbols of
\(\omega^{(n)}\) and \(\omega\) agree for all sufficiently
large \(n\).

For a finite word \(u=u_0\cdots u_{n-1}\), denote its length by
\(|u|=n\). We use \(uv\) for the concatenation of \(u\) and \(v\),
and use the same notation when \(v\) is infinite.
For \(\omega\in\Sig\) and \(n\geq1\), write
\(\omega|_n=\omega_0\cdots\omega_{n-1}\) for its prefix of length \(n\).
For \(1\leq N\leq n=|u|\), we denote the suffix of \(u\) of length \(N\) by
\[
 \suf_Nu:=u_{n-N}\cdots u_{n-1}.
\]
An infinite word \(\omega\) is called \emph{purely periodic} if
\(\omega=U^\infty\) for some nonempty finite word \(U\), that is,
\(\sigma^{|U|}\omega=\omega\). Similarly, it is called
\emph{eventually periodic} if some shift of it is purely periodic.

A nonempty subset \(\Omega\subset\Sig\) is called
a subshift if it is closed and \emph{shift invariant},
that is, \(\sigma(\Omega)\subseteq \Omega\).
For \(n\geq1\), we define its language of length \(n\) by
\(\Lang_n(\Omega)=\{\omega_j\cdots\omega_{j+n-1}:\omega\in \Omega,\ j\geq0\}\),
and write \(\Lang(\Omega)=\bigcup_{n\geq1}\Lang_n(\Omega)\).
Since \(\Omega\) is shift invariant, every word in \(\Lang_n(\Omega)\)
is the prefix of some sequence in \(\Omega\).
The cylinder and follower set of \(u\in\Lang(\Omega)\) are
\begin{equation}\label{eq:abstractcylinders}
 [u]_{\Omega}=\{\omega\in \Omega:\omega|_{|u|}=u\},\qquad
 \Fol_{\Omega}(u)=\{\eta\in\Sig:u\eta\in \Omega\}.
\end{equation}
In other words, \(\Fol_{\Omega}(u)\) consists of the infinite sequences
that can be added after \(u\) to obtain a sequence in \(\Omega\);
see Dokuchaev and Exel~\cite{DE}.
Each cylinder \([u]_{\Omega}\) is both open and closed
in \(\Omega\), and satisfies \(\operatorname{diam}[u]_{\Omega}\leq2^{-|u|}\).

\begin{definition}\label{def:SFT}
A subshift \(\Omega\) is said to be of \emph{finite type} if there
exists a finite set \(\mathcal B\) of finite words such that
\[
 \Omega=\{\omega\in\Sig:\text{no word of }\mathcal B
                                  \text{ occurs in }\omega\}.
\]
For convenience, we write SFT for a subshift of finite type.
The elements of \(\mathcal B\) are called forbidden words.
We call \(\Omega\) an \(N\)-step SFT if a sequence belongs
to \(\Omega\) exactly when all its blocks of length \(N+1\) belong to
\(\Lang_{N+1}(\Omega)\).
\end{definition}

\Needspace{7\baselineskip}
The following lemma is a one-sided version of the
characterization of SFTs in~\cite{LM}, stated in terms of follower sets.

\begin{lemma}[{\cite[Theorem~2.1.8]{LM}}]\label{lem:finitememory}
Let \(\Omega\subset\Sig\) be a subshift and let \(N\geq1\).
Then \(\Omega\) is an \(N\)-step SFT if and only if
\(\Fol_{\Omega}(u)=\Fol_{\Omega}(\suf_Nu)\) for every
\(u\in\Lang(\Omega)\) with \(|u|\geq N\).
\end{lemma}

We use \(\prec\) to denote the \emph{alternating lexicographic order} on
\(\Sig\). For two distinct sequences \(\omega,\nu\), let \(k\)
be the first index at which they differ.
Then \(\omega\prec\nu\) if and only if
\((-1)^k(\omega_k-\nu_k)<0\).
We write \(\omega\preceq\nu\) if \(\omega\prec\nu\) or
\(\omega=\nu\), and define \(\succ\) and \(\succeq\) similarly.
Throughout the paper, these symbols refer to alternating order.
Notice that deleting a common prefix preserves the order if its
length is even, and reverses the order if its length is odd.

\subsection{Kneading invariants}

Fix \((\beta,\alpha)\in\Delta_-\) and set \(c=\alpha/\beta\).
Recall that
\begin{equation*}
 T(x)=T_{-\beta,\alpha}(x)=
 \begin{cases}
  \alpha-\beta x,&0\leq x<c,\\
  \alpha+1-\beta x,&c<x\leq1,
 \end{cases}
\end{equation*}
thus \(T(c-)=0\), \(T(c+)=1\), \(T(0)=\alpha\),
and \(T(1)=1+\alpha-\beta\).
Denote by \(\pc_n\) the set of points whose orbits first reach
\(c\) at time \(n\). In particular, \(\pc_0=\{c\}\).
Set \(\pc=\bigcup_{n\geq0}\pc_n\).
Here the iterates are defined only until the first visit to \(c\).
Since each \(\pc_n\) is finite, \(\pc\) is countable, and hence
\((0,1)\setminus\pc\) is dense in \([0,1]\). For \(x\in(0,1)\)
with \(x\notin\pc\), we define the \emph{kneading sequence}
of such a point by
\[
 \tau(x)=\omega_0\omega_1\cdots,\qquad
 \omega_j=
 \begin{cases}1,&T^j(x)<c,\\2,&T^j(x)>c.\end{cases}
\]

\begin{definition}\label{def:naturalshift}
The \emph{lower kneading sequence} \(\tau(x-)\) and
the \emph{upper kneading sequence} \(\tau(x+)\) are defined by the
following one-sided limits, respectively, where \(y\in(0,1)\)
and \(y\notin\pc\):
\[
\left\{
\begin{aligned}
 \tau(x-)&:=\lim_{y\uparrow x,\;y\notin\pc}\tau(y),
                              &\qquad 0<x\leq1,\\[1ex]
 \tau(x+)&:=\lim_{y\downarrow x,\;y\notin\pc}\tau(y),
                              &\qquad 0\leq x<1.
\end{aligned}
\right.
\]
The limits are taken with respect to the metric \(d\).
The \emph{kneading invariants} are defined by
\((k_-,k_+):=(\tau(c-),\tau(c+))\).
 \end{definition}
 We use the one-sided convention of~\cite{DSS}. To see that the
limits exist, fix \(m\geq1\). The finite set
\(\bigcup_{j=0}^{m-1}\pc_j\) divides \((0,1)\) into intervals on
which the first \(m\) symbols are identical. Hence the limits in
Definition~\ref{def:naturalshift} exist. If \(x\in(0,1)\) and
\(x\notin\pc\), both limits equal \(\tau(x)\). If the orbit of
an endpoint avoids \(c\), its inward limit equals its branch coding.
If the orbit reaches \(c\), the inward limit still defines its
kneading sequence without specifying \(T(c)\).

Recall that the intermediate \((-\beta)\)-shift
of \(T_{-\beta,\alpha}\) is defined by
\begin{equation*}
 \Omega_{-\beta,\alpha}
 =\{\tau(x-):0<x\leq1\}\cup\{\tau(x+):0\leq x<1\}.
\end{equation*}
By definition, both kneading invariants belong to
\(\Omega_{-\beta,\alpha}\).
Since both branches reverse orientation, we have
 $\sigma (k_-)=\tau(0+),$ and
 $\sigma (k_+)=\tau(1-).$

The space \(\Omega_{-\beta,\alpha}\) is a subshift. Indeed, let
\(\omega^{(m)}\in\Omega_{-\beta,\alpha}\) converge to \(\omega\in\Sig\).
Choose \(x_m\in(0,1)\setminus\pc\) such that \(\tau(x_m)\) and
\(\omega^{(m)}\) have the same first \(m\) symbols.
By compactness of \([0,1]\), we may pass to a subsequence such that
\(x_m\) approaches some \(x\) from the left or from the right, or
\(x_m=x\) for all \(m\). Then \(\omega\) equals \(\tau(x-)\),
\(\tau(x+)\), or \(\tau(x)\), respectively. In each case,
\(\omega\in\Omega_{-\beta,\alpha}\), so the space is closed.
Since each branch is decreasing, shifting a one-sided coding gives
a one-sided coding at its branch image, with the side reversed.
Hence \(\sigma(\Omega_{-\beta,\alpha})\subseteq\Omega_{-\beta,\alpha}\).

\begin{proposition}\label{prop:kneadingdescription}
For the intermediate \((-\beta)\)-shift
\(\Omega_{-\beta,\alpha}\), we have
\begin{equation}\label{eq:kneadingdescription}
 \Omega_{-\beta,\alpha}
 =\left\{\omega\in\Sig:
 \sigma(k_-)\preceq\sigma^n\omega\preceq\sigma(k_+),\quad
 \forall n\geq0\right\}.
\end{equation}
\end{proposition}

\begin{proof}
Let \(x<y\) belong to \((0,1)\setminus\pc\).
Since \(T\) is expanding on each branch, there is a least
\(j\geq0\) such that \(T^j(x)\) and \(T^j(y)\) lie on opposite
sides of \(c\). The order of \(T^j(x),T^j(y)\)
is reversed exactly when \(j\) is odd, so \(\tau(x)\prec\tau(y)\).
Taking one-sided limits and using the same expansion argument, the codings at
distinct points remain strictly ordered. Hence the minimum and maximum
codings are \(\sigma(k_-)=\tau(0+)\) and \(\sigma(k_+)=\tau(1-)\),
respectively. By shift invariance, every
sequence in \(\Omega_{-\beta,\alpha}\) satisfies the conditions
in~\eqref{eq:kneadingdescription}.

Conversely, suppose that \(\omega\) satisfies the
shift conditions but \(\omega\notin\Omega_{-\beta,\alpha}\).
By the condition at \(n=0\), \(\omega\) lies between the minimum
and maximum codings. Since the space is compact and inequalities in
alternating order are preserved under limits, there are a largest
coding \(a\prec\omega\) and a smallest coding
\(b\succ\omega\). They correspond to the same point; otherwise,
a point outside \(\pc\) strictly between the corresponding points
would have its coding between \(a\) and \(b\), a contradiction.
Since endpoints and points outside \(\pc\) have only one coding,
this common point belongs to \(\pc_n\cap(0,1)\) for some \(n\geq0\).
Let \(w\) be its word of length \(n\) before reaching \(c\).
Then \(\{a,b\}=\{wk_-,wk_+\}\), and \(\omega\) also has prefix
\(w\). Deleting this prefix, with order reversed when \(n\) is odd,
gives \(k_-\prec\sigma^n\omega\prec k_+\).
Since \(k_-=1\sigma(k_-)\) and \(k_+=2\sigma(k_+)\), deleting
the first symbol once more gives
\(\sigma^{n+1}\omega\prec\sigma(k_-)\) if \(\omega_n=1\), and
\(\sigma^{n+1}\omega\succ\sigma(k_+)\) if \(\omega_n=2\).
Both cases contradict the shift condition at \(n+1\).
\end{proof}

\subsection{Kneading determinant}
Define the address map on \( \Omega_{-\beta,\alpha}\) by
\begin{equation}\label{eq:address}
 \pi(\omega)=\frac{\alpha}{\beta+1}
       +\sum_{j=0}^{\infty}\frac{(-1)^j(\omega_j-1)}{\beta^{j+1}}.
\end{equation}
Since the series converges uniformly, \(\pi\) is continuous.
By shifting the series, for every \(\omega\in\Omega_{-\beta,\alpha}\), we have
\begin{equation}\label{eq:branchidentity}
 \pi(\sigma\omega)=\alpha-\beta\pi(\omega)+\omega_0-1.
\end{equation}
For \(x\in(0,1)\), \(x\notin\pc\), write
\(\omega=\tau(x)\) and \(x_j=T^j(x)\). Iterating
\(x_j=(\alpha+\omega_j-1-x_{j+1})/\beta\) yields
\[x=\sum_{j=0}^{n-1}\frac{(-1)^j(\alpha+\omega_j-1)}{\beta^{j+1}}
       +\frac{(-1)^n x_n}{\beta^n}.\]
Letting \(n\to\infty\), we obtain \(\pi(\tau(x))=x\).
By continuity, we also have \(\pi(\tau(x-))=x\) and
\(\pi(\tau(x+))=x\) whenever these one-sided sequences are defined.
Hence \(\pi(\Omega_{-\beta,\alpha})=[0,1]\) and
\(\pi(k_-)=\pi(k_+)=c\).

For the kneading invariants \((k_-,k_+)\) of \(T_{-\beta,\alpha}\),
define the kneading determinant by
\[
 K_{(k_-,k_+)}(z)=\sum_{j=0}^{\infty}
       \bigl((k_+)_j-(k_-)_j\bigr)z^j,
 \qquad |z|<1.
\]
The coefficients are bounded, so this function is analytic in
the open unit disk. Since \((k_-)_0=1\) and \((k_+)_0=2\),
we have \(K_{(k_-,k_+)}(0)=1\). By~\eqref{eq:address} and
\(\pi(k_-)=\pi(k_+)=c\), we obtain
\begin{equation}\label{eq:kneadingzero}
 K_{(k_-,k_+)}(-1/\beta)=0.
\end{equation}
Identity~\eqref{eq:kneadingzero} will be used in
Section~\ref{sec:constructive}.

\section{Proof of Theorem \ref{main:criterion}}\label{sec:criterion}

We first consider the intervals corresponding to symbolic
cylinders for \(T_{-\beta,\alpha}\).
For \(\Omega_{-\beta,\alpha}\), we write
\([u]=[u]_{\Omega_{-\beta,\alpha}}\) and
\(\Fol(u)=\Fol_{\Omega_{-\beta,\alpha}}(u)\),
as in~\eqref{eq:abstractcylinders}.
For \(i\in\{1,2\}\), let
\(f_i(x)=\alpha+i-1-\beta x\), and set \(I_1=(0,c)\), \(I_2=(c,1)\).
For a finite word \(u=u_0\cdots u_{\ell-1}\), define its open geometric
cylinder by
\[
 J(u)=\{x\in(0,1):T^j(x)\in I_{u_j},\ 0\leq j<\ell\}.
\]
The iterates are determined by these branch inequalities
without using \(T(c)\). Let
\[
 \Phi_u=f_{u_{\ell-1}}\circ\cdots\circ f_{u_0}.
\]
Then \(\Phi_u\) is an affine map with slope
\((-\beta)^\ell\), and it agrees with \(T^\ell\) on \(J(u)\).

\begin{lemma}\label{lem:cylinders}
A finite word \(u\) is allowed in \(\Omega_{-\beta,\alpha}\)
if and only if \(J(u)\ne\varnothing\). In this case, \(J(u)\)
is an open interval of diameter at most \(\beta^{-|u|}\), and
\[
 \pi([u])=\cl{J(u)},\qquad
 \pi(\Fol(u))=\Phi_u(\cl{J(u)}).
\]
\end{lemma}

\begin{proof}
Suppose first that \(J(u)\ne\varnothing\).
Since the branch conditions defining \(J(u)\) are strict affine
inequalities, \(J(u)\) is an open interval. Hence we can choose
a point in \(J(u)\setminus\pc\), whose coding begins with \(u\).
Thus \(u\) is allowed. Conversely, if \(u\)
is allowed, then some one-sided coding begins with \(u\).
By the definition of the one-sided limit, a sufficiently close
point \(y\in(0,1)\setminus\pc\) on that side has the same
prefix. Thus \(y\in J(u)\), proving the equivalence.
By the same argument, every sequence in \([u]\) can be
approximated by codings of points in \(J(u)\setminus\pc\).
Since \([u]\) is closed, we obtain
\begin{equation*}
 [u]=\cl{\{\tau(x):x\in J(u),\ x\notin\pc\}}.
\end{equation*}
By the continuity of \(\pi\), compactness, and the identity
\(\pi(\tau(x))=x\), we obtain \(\pi([u])=\cl{J(u)}\).
Let \(\ell=|u|\). Since \(\Fol(u)=\sigma^\ell[u]\), iterating
the identity \(\pi(\sigma\omega)=\alpha-\beta\pi(\omega)+\omega_0-1\)
gives
\[
 \pi(\Fol(u))=\pi(\sigma^\ell[u])
 =\Phi_u(\pi([u]))=\Phi_u(\cl{J(u)}).
\]
Finally, since \(\Phi_u(J(u))\subset(0,1)\) and the
modulus of its slope is \(\beta^\ell\), we obtain the required
diameter bound.
\end{proof}

Using Lemmas~\ref{lem:finitememory} and~\ref{lem:cylinders},
we now prove the equivalence of conditions~\((1)\)
and~\((2)\) in Theorem~\ref{main:criterion}.
\begin{proposition}\label{prop:endpointreturn}
Fix $1<\beta<2, \ \beta-1<\alpha<1,$  and $c=\frac{\alpha}{\beta}.$
Then \(\Omega_{-\beta,\alpha}\) is an SFT if and only if both
endpoints reach \(c\).
\end{proposition}

\begin{proof}
Suppose first that
\(\Omega_{-\beta,\alpha}\) is an \(N\)-step SFT with \(N\geq1\).
By Lemma~\ref{lem:finitememory}, for every allowed word \(u\)
of length at least \(N\), we have
\begin{equation}\label{eq:finitefol}
 \Fol(u)=\Fol(\suf_Nu).
\end{equation}
We prove that both endpoints reach \(c\) by contradiction.
Suppose that \(e\in\{0,1\}\) never reaches \(c\).
Since \((\beta,\alpha)\in\Delta_-\), we have
\(T(0)=\alpha\in(0,1)\) and
\(T(1)=1+\alpha-\beta\in(0,1)\). Neither branch takes the value
\(0\) or \(1\) away from the discontinuity. Hence the positive orbit
of \(e\) lies in \((0,1)\). It follows that every finite prefix of the
endpoint kneading sequence is realized by nearby points \(x\in(0,1)\)
with \(x\notin\pc\).

Let \(u\) be its prefix of length \(N+1\), let
\(v=\suf_Nu\), and set \(z=T^{N+1}(e)\). Since \(J(u)\) is adjacent
to \(e\), the point \(e\) is an endpoint of \(\cl{J(u)}\).
Therefore, the nondegenerate closed interval
\(\pi(\Fol(u))=\Phi_u(\cl{J(u)})\)
has \(z\) as an endpoint. On the other hand, the next
\(N\) branch inequalities for \(T(e)\) are strict, so \(T(e)\) lies
in the interior of \(J(v)\). Hence
\(z=\Phi_v(T(e))\) is an interior point of
\(\pi(\Fol(v))=\Phi_v(\cl{J(v)})\).
Thus these two intervals are different, which contradicts
\eqref{eq:finitefol}. Applying the same argument to each endpoint,
we conclude that both endpoints reach \(c\).

Conversely, suppose both endpoints reach \(c\), and let
\(n_0,n_1\) be their first hitting times.
Set \(M=\max(n_0,n_1)\), let \(r=M+1\), and define
the finite set \(B=\{0,1\}\cup\bigcup_{j=0}^{M}\pc_j\).
If \(x\) first reaches \(c\) at time \(j\geq1\),
then \(T(x)\) first reaches \(c\) at time \(j-1\).
The images of \(0,1\) also belong to \(B\), and
\(T(c-)=0\), \(T(c+)=1\).
Hence \(B\) is forward invariant when both branch limits are included.

We next show that the nonempty cylinders \(J(u)\),
\(|u|=r\), are exactly the components of \((0,1)\setminus B\).
The symbols through time \(r-1=M\) are constant on each component.
Moreover, each point of \(B\cap(0,1)\) separates
distinct prefixes. Combining this with the fact that each prefix cylinder is
an interval by Lemma~\ref{lem:cylinders}, we obtain the assertion.

For each partition interval \(J(u)\), the map \(T\) is monotone and
the endpoints of \(T(J(u))\) belong to \(B\). Thus \(T(J(u))\) is
a union of partition intervals and finitely many points of
\(B\). We also have
\(T(J(u)\setminus\pc)=T(J(u))\setminus\pc\).
Indeed, for \(x\in J(u)\), we have \(x\notin\pc\)
if and only if \(T(x)\notin\pc\). Since both endpoints reach \(c\),
we have \(B\subseteq\pc\). Thus \(T(J(u))\setminus\pc\) is the union of
the sets \(J(v)\setminus\pc\) over the partition intervals \(J(v)\)
contained in \(T(J(u))\). For \(x\in J(u)\setminus\pc\), we have
\(\tau(T(x))=\sigma\tau(x)\). Taking closures in the sequence space
as in the proof of Lemma~\ref{lem:cylinders}, and using compactness
and the continuity of \(\sigma\), we obtain that \(\sigma[u]\) is the
union of the corresponding length-\(r\) cylinders \([v]\).

Let \(\omega\in\Sig\) have all its length-\(r+1\)
blocks in \(\Lang_{r+1}(\Omega_{-\beta,\alpha})\), and put
\(u_j=\omega_j\cdots\omega_{j+r-1}\). Since
\(\omega_j\cdots\omega_{j+r}\) is allowed, we have
\(\sigma[u_j]\cap[u_{j+1}]\ne\varnothing\).
By the preceding observation and the disjointness of length-\(r\) cylinders,
\([u_{j+1}]\subseteq\sigma[u_j]\).
For any \(m\geq1\), choose a sequence in \([u_m]\).
By the above inclusions, we can successively choose its preimages
under \(\sigma\) in \([u_{m-1}],\ldots,[u_0]\).
The resulting sequence in \([u_0]\) has the same first \(m+r\)
symbols as \(\omega\). Since \(m\) can be arbitrarily large and
\(\Omega_{-\beta,\alpha}\) is closed, we have
\(\omega\in\Omega_{-\beta,\alpha}\).
Thus the shift is an \(r\)-step SFT.
\end{proof}
Next we express the condition that both endpoints reach
\(c\) in terms of the kneading invariants.

\begin{proposition}\label{prop:kneadingreturn}
Both endpoints of \(T_{-\beta,\alpha}\) reach \(c\)
if and only if there exist \(p,q\geq1\) such that
\(\sigma^p k_-\in\{k_-,k_+\}\) and
\(\sigma^q k_+\in\{k_-,k_+\}\).
\end{proposition}

\begin{proof}
Suppose both endpoints reach \(c\), and let
\(n_0,n_1\) be their first hitting times. Put \(p=n_0+1\) and
\(q=n_1+1\). Since the one-sided images of \(c\) are \(0+\)
and \(1-\), the side-reversal argument in
\cite[Proposition~3.2]{DSS} gives
\begin{equation}\label{eq:returnparity}
 \sigma^pk_-=
 \begin{cases}k_-,&p\text{ even},\\k_+,&p\text{ odd},\end{cases}
 \qquad
 \sigma^qk_+=
 \begin{cases}k_+,&q\text{ even},\\k_-,&q\text{ odd}.\end{cases}
\end{equation}

Conversely, suppose that positive shifts of both
kneading invariants belong to \(\{k_-,k_+\}\).
By \(\pi(k_-)=\pi(k_+)=c\) and~\eqref{eq:branchidentity}, each
such shift corresponds to an orbit reaching \(c\).
After the initial step from \(c\), this orbit follows that of
the corresponding endpoint until it first reaches \(c\).
Hence both endpoints reach \(c\), and the first shifts returning
to \(\{k_-,k_+\}\) occur at \(n_0+1\) and \(n_1+1\), respectively.
\end{proof}

\begin{proof}[Proof of Theorem~\ref{main:criterion}]
By Proposition~\ref{prop:endpointreturn}, conditions~\((1)\)
and~\((2)\) are equivalent. By Proposition~\ref{prop:kneadingreturn},
conditions~\((2)\) and~\((3)\) are equivalent. Combining these two
propositions, we complete the proof.
\end{proof}

\begin{remark}\label{rem:fullshiftboundary}
At the boundary parameter pair \((\beta,\alpha)=(2,1)\), we have
\(\Omega_{-2,1}=\Sig\), which is an SFT. However, the endpoints
form the cycle \(0\mapsto1\mapsto0\) and never reach \(c=1/2\).
Thus the return criterion in Theorem~\ref{main:criterion}
does not extend to this parameter pair.
\end{remark}

Suppose both endpoints reach \(c\). Let \(p,q\)
be the least positive integers such that
\(\sigma^p k_-\in\{k_-,k_+\}\) and
\(\sigma^q k_+\in\{k_-,k_+\}\), respectively.
Denote the corresponding \emph{first-return words}
by \(U=(k_-)|_p\) and \(V=(k_+)|_q\).
By~\eqref{eq:returnparity}, the kneading invariants are given in
Table~\ref{tab:returnparity}.
For positive slopes, an intermediate \(\beta\)-shift
is an SFT if and only if both kneading invariants are purely
periodic~\cite[Theorem~1.3]{LSS}. For negative slopes, this condition
is sufficient but not necessary; see Example~\ref{ex:mixed}.
\begin{table}[!htbp]
\centering
\caption{Kneading invariants for the four parity cases.}
\label{tab:returnparity}
\begin{tabular}{cccc}
\toprule
\(p\) & \(q\) & \(k_-\) & \(k_+\)\\
\midrule
even & even & \(U^\infty\) & \(V^\infty\)\\
odd  & even & \(UV^\infty\) & \(V^\infty\)\\
even & odd  & \(U^\infty\) & \(VU^\infty\)\\
odd  & odd  & \((UV)^\infty\) & \((VU)^\infty\)\\
\bottomrule
\end{tabular}
\end{table}

\section{Proof of Theorem \ref{main:density}}\label{sec:constructive}

We first find parameters for which the endpoint \(0\) reaches
\(c\). Next we vary the parameters along a curve on which this return
is unchanged, and prove that parameters for which the endpoint \(1\)
also reaches \(c\) are dense on this curve.
Throughout this section, write \(k_\pm(b,a)\) for the kneading
invariants of \(T_{-b,a}\).

\subsection{Endpoint returns}

For \(b>1\) and \(n\geq1\), put
\(S_n(b)=\sum_{j=0}^{n-1}(-b)^j=(1-(-b)^n)/(1+b)\).

\begin{lemma}\label{lem:singlereturndense}
Fix \(\beta\in(1,2)\). The set of
\(\alpha\in(\beta-1,1)\) for which the endpoint \(0\) reaches
\(c_{\beta,\alpha}=\alpha/\beta\) is dense in \((\beta-1,1)\).
\end{lemma}

\begin{proof}
We prove this by contradiction. Suppose that there is
a nonempty open interval \(J\subset(\beta-1,1)\) such that,
for every \(\alpha\in J\), the orbit of \(0\) under
\(T_{-\beta,\alpha}\) never reaches the discontinuity \(c_{\beta,\alpha}\).
We show by induction that the branch symbol at each
time is the same for all \(\alpha\in J\).
At time \(0\), the symbol is \(1\) since \(0<c_{\beta,\alpha}\).
Suppose that the symbols through
time \(j-1\) are fixed on \(J\), where \(j\geq1\). Composing these
branches makes \(T_{-\beta,\alpha}^j(0)\) an affine function of
\(\alpha\). Hence \(T_{-\beta,\alpha}^j(0)-\alpha/\beta\) is
continuous on \(J\) and never vanishes by our assumption.
By the intermediate value theorem, it has the same sign throughout
\(J\), so the branch symbol at time \(j\) is also fixed.
This proves that every finite itinerary is constant on \(J\).

Write \(w_0w_1\cdots\) for this common branch coding.
For every \(n\geq1\), composing the corresponding branches gives
\[
 T_{-\beta,\alpha}^n(0)=\alpha S_n(\beta)+D_n,
 \qquad \alpha\in J,
\]
where
\(D_n=\sum_{j=0}^{n-1}(w_j-1)(-\beta)^{n-1-j}\)
is independent of \(\alpha\). Choose distinct
\(\alpha_1,\alpha_2\in J\). Both orbits remain in \([0,1]\), so
\[
 |S_n(\beta)|\,|\alpha_1-\alpha_2|\leq1
 \qquad\text{for every }n\geq1.
\]
This contradicts \(|S_n(\beta)|\to\infty\).
\end{proof}

\begin{remark}\label{rem:onereturndense}
The same argument shows that, for each fixed
\(\beta\in(1,2)\), the set of \(\alpha\in(\beta-1,1)\) for which
the endpoint \(1\) reaches \(c_{\beta,\alpha}\) is also dense in
\((\beta-1,1)\).
\end{remark}

Suppose that the endpoint \(0\) first reaches
\(c_{\beta_0,\alpha_0}\) at time \(n\geq1\) for
\((\beta_0,\alpha_0)\in\Delta_-\). Let \(w=w_0\cdots w_{n-1}\)
be its branch word before this visit.
The branch inequalities at times \(0,\ldots,n-1\) are strict.
By continuity, these finitely many inequalities remain valid in a
sufficiently small neighbourhood of \((\beta_0,\alpha_0)\) contained
in \(\Delta_-\). For parameters \((b,a)\) in this neighbourhood,
the endpoint \(0\) therefore follows the same word \(w\), and
composing the corresponding branches gives
\(T_{-b,a}^n(0)=aS_n(b)+D_w(b)\), where
\(D_w(b)=\sum_{j=0}^{n-1}(w_j-1)(-b)^{n-1-j}\).
Thus the condition \(T_{-b,a}^n(0)=c_{b,a}=a/b\) is equivalent
to \(aS_n(b)+D_w(b)=a/b\). Solving for \(a\), we obtain
\(a=\varphi(b)\), where
\[
 \varphi(b)=\frac{D_w(b)}{1/b-S_n(b)}.
\]
Since \(1/b-S_n(b)=(1+b(-b)^n)/(b(b+1))\ne0\) for \(b>1\),
the function \(\varphi\) is real analytic. The original parameter
satisfies the return equation, so \(\varphi(\beta_0)=\alpha_0\).
Choose an open interval \(I\) containing \(\beta_0\) small enough
that \((b,\varphi(b))\) stays in the above neighbourhood for every
\(b\in I\). Then \(0\) follows the same word \(w\) and first reaches
\(c_{b,\varphi(b)}=\varphi(b)/b\) at time \(n\).

\Needspace{6\baselineskip}
\begin{lemma}\label{lem:returncurve}
The values of \(b\in I\) for which the endpoint
\(1\) reaches \(c_{b,\varphi(b)}\) under \(T_{-b,\varphi(b)}\)
form a dense subset of \(I\).
\end{lemma}

\begin{proof}
Let \(U=1w\) and \(p=n+1\). Applying the same argument about
the reversal of sides as in the proof of
Proposition~\ref{prop:kneadingreturn} to the endpoint \(0\), we have
\begin{equation}\label{eq:singleparity}
 k_-(b,\varphi(b))=
 \begin{cases}
 U^\infty,&p\text{ even},\\
 U k_+(b,\varphi(b)),&p\text{ odd}.
 \end{cases}
\end{equation}
We prove the density assertion by contradiction.
Suppose that there is a nonempty open subinterval \(J\subset I\) such that the endpoint
\(1\) never reaches \(c_{b,\varphi(b)}\) for any \(b\in J\).
As in the proof of Lemma~\ref{lem:singlereturndense}, the absence
of a return to the discontinuity implies that each finite itinerary
of \(1\) is locally constant on \(J\). Since \(J\) is an interval,
each such itinerary is constant on \(J\).
Thus \(k_+(b,\varphi(b))\) is constant.
By~\eqref{eq:singleparity}, \(k_-(b,\varphi(b))\) is also constant.
Hence their kneading determinant is a fixed function \(K\),
analytic on \(|z|<1\), with \(K(0)=1\).
By~\eqref{eq:kneadingzero}, we have \(K(-1/b)=0\) for every \(b\in J\).
These zeros have an accumulation point in the open unit disk.
The identity theorem therefore implies that \(K\) vanishes identically there,
contradicting \(K(0)=1\). This proves the density assertion.
\end{proof}

By Remark~\ref{rem:onereturndense}, we may also start with
a return of the endpoint \(1\). The same construction and proof apply with the roles of
\(0\) and \(1\) interchanged.

\subsection{Simultaneous approximation}

We now approximate both the parameters and their kneading invariants.

\begin{proposition}\label{prop:kneadingapproximation}
For every \((\beta,\alpha)\in\Delta_-\), \(N\geq1\) and
\(\varepsilon>0\), there exists \((b,a)\in\mathcal F_-\) such that
\(|b-\beta|+|a-\alpha|<\varepsilon\),
\(k_-(b,a)|_N=k_-(\beta,\alpha)|_N\), and
\(k_+(b,a)|_N=k_+(\beta,\alpha)|_N\).
\end{proposition}

\begin{proof}
If both endpoints reach \(c\), then we can take
\((b,a)=(\beta,\alpha)\) by Theorem~\ref{main:criterion}.
It remains to consider the following two cases.

\emph{Case 1. Neither endpoint reaches \(c\).}
The finite endpoint orbit segments determining the first \(N\) symbols
of both kneading invariants satisfy strict branch inequalities.
By continuity of the branches and of \(c=\alpha/\beta\), there is
a neighborhood \(W\subset\Delta_-\) of \((\beta,\alpha)\) in which
these prefixes are unchanged and the required parameter bound holds.
By Lemma~\ref{lem:singlereturndense}, choose
\((\beta,\alpha_0)\in W\) for which \(0\) reaches \(c\).
If \(1\) also reaches \(c\), this parameter has the required
properties. Otherwise, consider the curve constructed above
through \((\beta,\alpha_0)\) and restrict it to \(W\).
By Lemma~\ref{lem:returncurve}, choose a parameter on this curve
for which \(1\) also reaches \(c\).
Theorem~\ref{main:criterion} then gives an SFT parameter in \(W\).

\emph{Case 2. Exactly one endpoint reaches \(c\).}
We only consider the case that \(0\) reaches \(c\), since the
other case can be proved similarly.
Consider the above curve through \((\beta,\alpha)\),
and let \(U\) and \(p\) be the fixed word and length
in~\eqref{eq:singleparity}.
Since \(1\) does not reach \(c\) at the original parameter,
the first \(N\) symbols of \(k_+\) are unchanged at all
sufficiently close parameters.
By~\eqref{eq:singleparity}, the first \(N\) symbols
of \(k_-\) are also unchanged at sufficiently close parameters
on the return curve.
Indeed, if \(p\) is even, then \(k_-=U^\infty\) is unchanged.
If \(p\) is odd, then \(k_-\) is obtained by attaching the fixed
prefix \(U\) to \(k_+\).
Choose a sufficiently short curve segment so that both prefixes
and the required parameter bound hold. By Lemma~\ref{lem:returncurve},
this segment contains a parameter for which \(1\) also reaches
\(c\). Theorem~\ref{main:criterion} completes the proof.
\end{proof}

\begin{proof}[Proof of Theorem~\ref{main:density}]
Fix \((\beta,\alpha)\in\Delta_-\). For each \(i\geq1\),
apply Proposition~\ref{prop:kneadingapproximation}
with \(N=i\) and \(\varepsilon=1/i\).
The resulting SFT parameters \((b_i,a_i)\) converge to
\((\beta,\alpha)\), and their kneading invariants converge to
\(k_-(\beta,\alpha)\) and \(k_+(\beta,\alpha)\), respectively.
In particular, \(\mathcal F_-\) is dense in \(\Delta_-\).
\end{proof}

We conclude with an example illustrating that an intermediate
\((-\beta)\)-shift can be of finite type without both kneading
invariants being purely periodic.

\begin{example}\label{ex:mixed}
Let \(\beta\in(1,2)\) be the unique root of
\(\beta^4-\beta^3-1=0\), and set
\(\alpha=\beta/(\beta^2-\beta+1)\). Then
\((\beta,\alpha)\in\Delta_-\). A direct calculation shows that
the endpoints \(0\) and \(1\) first reach \(c\) after two and three
iterates, respectively. The corresponding first-return words are
\(U=112\) and \(V=2211\). By Table~\ref{tab:returnparity},
\[
 k_-=112(2211)^\infty,\qquad k_+=(2211)^\infty.
\]
By Proposition~\ref{prop:endpointreturn}, the shift is an
SFT. However, the block \(222\) occurs exactly once in \(k_-\),
so \(k_-\) is not purely periodic.
\end{example}

\medskip
\noindent\textbf{Acknowledgements.}
The authors thank Qian Xiao for helpful discussions.
Yiming Ding was supported by the National Natural Science Foundation
of China (No.~12271418).
Yun Sun was supported by the National Natural Science
Foundation of China (No.~12601157).

\medskip
\noindent\textbf{Statements and AI declarations.}
The authors declare no conflict of interest.
No datasets were generated or analyzed in this study.


\begin{thebibliography}{99}
\urlstyle{same}

\bibitem{BSV}
M.~Barnsley, W.~Steiner, and A.~Vince.
Critical itineraries of maps with constant slope and one discontinuity.
\emph{Math. Proc. Cambridge Philos. Soc.}, 157(3):547--565, 2014.

\bibitem{Blanchard}
F.~Blanchard.
\(\beta\)-expansions and symbolic dynamics.
\emph{Theoret. Comput. Sci.}, 65(2):131--141, 1989.

\bibitem{DR}
K.~Dajani and S.~D.~Ramawadh.
Symbolic dynamics of \((-\beta)\)-expansions.
\emph{J. Integer Seq.}, 15(2):Article~12.2.6, 2012.

\bibitem{DSS}
Y.~Ding, W.~Steiner, and Y.~Sun.
Open dynamical systems in topologically expansive decreasing Lorenz maps.
Preprint, \href{https://arxiv.org/abs/2610.03327}{arXiv:2610.03327}, 2026.

\bibitem{DE}
M.~Dokuchaev and R.~Exel.
Partial actions and subshifts.
\emph{J. Funct. Anal.}, 272(12):5038--5106, 2017.

\bibitem{FL}
C.~Frougny and A.~C.~Lai.
Negative bases and automata.
\emph{Discrete Math. Theor. Comput. Sci.}, 13(1):75--94, 2011.

\bibitem{HLY}
H.~Hu, Z.~Li, and Y.~Yu.
A note on \((-\beta)\)-shifts with the specification property.
\emph{Publ. Math. Debrecen}, 91(1--2):123--131, 2017.

\bibitem{HLS}
Y.~Huang, L.~Liao, and Y.~Sun.
Negative \(\beta\)-transformations: invariant measures, matching
properties and associated subshifts.
Preprint, \href{https://arxiv.org/abs/2603.14862}{arXiv:2603.14862}, 2026.

\bibitem{HS}
J.~H.~Hubbard and C.~T.~Sparrow.
The classification of topologically expansive Lorenz maps.
\emph{Comm. Pure Appl. Math.}, 43(4):431--443, 1990.

\bibitem{IS}
S.~Ito and T.~Sadahiro.
Beta-expansions with negative bases.
\emph{Integers}, 9(3):239--259, 2009.

\bibitem{LSS}
B.~Li, T.~Sahlsten, and T.~Samuel.
Intermediate beta-shifts of finite type.
\emph{Discrete Contin. Dyn. Syst.}, 36(1):323--344, 2016.

\bibitem{LSSS}
B.~Li, T.~Sahlsten, T.~Samuel, and W.~Steiner.
Denseness of intermediate \(\beta\)-shifts of finite-type.
\emph{Proc. Amer. Math. Soc.}, 147(5):2045--2055, 2019.

\bibitem{LS}
L.~Liao and W.~Steiner.
Dynamical properties of the negative beta-transformation.
\emph{Ergodic Theory Dynam. Systems}, 32(5):1673--1690, 2012.

\bibitem{LM}
D.~Lind and B.~Marcus.
\emph{An Introduction to Symbolic Dynamics and Coding}.
Cambridge University Press, Cambridge, 1995.

\bibitem{NNZ}
F.~Nguema Ndong.
Zeta function and negative beta-shifts.
\emph{Monatsh. Math.}, 188(4):717--751, 2019.

\bibitem{NNBM}
F.~Nguema Ndong and A.~Bertrand-Mathis.
Exchange of intervals and intrinsic ergodicity of the negative beta shift.
Preprint, \href{https://arxiv.org/abs/2503.06597v5}{arXiv:2503.06597v5}, 2026.

\bibitem{Parry}
W.~Parry.
On the \(\beta\)-expansions of real numbers.
\emph{Acta Math. Acad. Sci. Hungar.}, 11(3--4):401--416, 1960.

\bibitem{Renyi}
A.~R\'enyi.
Representations for real numbers and their ergodic properties.
\emph{Acta Math. Acad. Sci. Hungar.}, 8(3--4):477--493, 1957.

\bibitem{SY}
M.~Shinoda and K.~Yamamoto.
Intrinsic ergodicity for factors of \((-\beta)\)-shifts.
\emph{Nonlinearity}, 33(1):598--609, 2020.

\bibitem{Steiner}
W.~Steiner.
Digital expansions with negative real bases.
\emph{Acta Math. Hungar.}, 139(1--2):106--119, 2013.

\bibitem{SLD}
Y.~Sun, B.~Li, and Y.~Ding.
Fiber denseness of intermediate \(\beta\)-shifts of finite type.
\emph{Nonlinearity}, 36(11):5973--5997, 2023.

\end{thebibliography}
\end{document}